\documentclass[smallextended,envcountsect]{svjour3}
\smartqed
\usepackage{graphicx}
\usepackage{amsfonts}
\usepackage{amsmath}
\usepackage[usenames,dvipsnames]{pstricks}
\usepackage{epsfig}
\usepackage{bbding}

\newcommand{\nid}{\noindent}
\newcommand{\bc}{\begin{center}}
\newcommand{\br}{\begin{right}}
\newcommand{\ec}{\end{center}}
\newcommand{\be}{\begin{equation}}
\newcommand{\ee}{\end{equation}}

\newcommand{\vl}{\mid}

\newcommand{\rar}{\rightarrow}
\newcommand{\grad}{\nabla}

\newcommand{\p}{\partial}

\newcommand{\bint}{\mbox{int} \,}

\newcommand{\bco}{\mbox{conv} \,}

\newcommand{\bmin}{\mbox{min} \,}

\newcommand{\exs}{\exists}

\newcommand{\meq}{\geq}

\newcommand{\al}{\alpha}

\newcommand{\eps}{\varepsilon}

\newtheorem{thm}{Theorem}[section]

\newtheorem{cor}{Corollary}[section]

\newtheorem{defi}{Definition}[section]
\newtheorem{ex}{Example}[section]

\journalname{JOTA}

\begin{document}

\title{Calculation of subdifferentials and
codifferentials. Construction of continuous codifferentials. I}


\titlerunning{Calculation of subdifferentials and
codifferentials}
\dedication{Dedicated to the memory of my teacher Prof. V.F.
Demyanov, who invented codifferentials}

\author{Igor M. Prudnikov}

\institute{Igor Mihailovich Prudnikov \at
             Scientific Center of Smolensk Federal Medical  University, Smolensk, Russia, 214000\\
              \email{pim\underline{ }10@hotmail.com}
}

\date{Received: date / Accepted: date}

\maketitle

\begin{abstract}
The author studies the twice codifferentiable functions introduced
by Professor V.F. Demyanov, and how to calculate their
subdifferentials and codifferentials.  The simpler case, when a
function is twice hypo-differentiable, is initially considered.
There is proved that twice hypo-differentiable positively
homogeneous function of the second order is the maximum of
quadratic forms over some set of matrices that coincides with the
convex hull of the limit matrices, that are calculated at the
points where the original function is twice differentiable, and
the points tend to zero. There is shown that the set of limit
matrices coincides with the subdifferential of the second order,
introduced by the author, of the original functions at zero. The
first and second subdifferentials are used to compute the second
codifferential at a point. Moreover, the second hypodifferential
and hyperdifferential are calculated up to equivalence. The proved
theorems give the rules for calculation of the subdifferentials
and continuous codifferentials. It is important for practical
optimization.
\end{abstract}



\section{Introduction}

Many authors introduce generalized gradients and matrices for
Lipschitz functions in different ways \cite{clark1} -
\cite{morduchrockafel}. But unlike the smooth case, the
subdifferentials, consisting of the generalized gradients, are not
continuous in the Hausdorff metric. Therefore, pairs of some sets,
that are analogues of the subdifferentials, were introduced in
\cite{demrub}, which are continuous in the Hausdorff metric for
wide class of functions. We will study the Demyanov-Rubinov
subdifferentials and their connection with the codifferentials.

Demyanov V.F. and Rubinov A.M. have introduced \cite{demrub}
codifferentiable and twice codifferentiable functions. They called
$ f (\cdot): \mathbb{R}^ n \rar \mathbb{R} $ {\em
codifferentiable} at a point $ x $ if there exist convex compact
sets $ \underline {d} f (x), \overline {d} f (x) $ from $
\mathbb{R}^ {n + 1} $, called the {\em hypodifferential} and {\em
hyperdifferential} respectively, for which the decomposition  \be
f (x + \Delta) = f (x) + \max_ {[a, v] \in \underline {d} f (x)}
[a + (v, \Delta)] + \min_ {[b, w] \in \overline {d} f (x)} [b +
(w, \Delta)] + o_x (\Delta), \label {quasidiff1} \ee is true and
{\em twice codifferentiable} at a point $ x $ if there exist
convex compacts $ \underline {d} ^ 2 f (x), \overline {d} ^ 2 f
(x) $ of $ \mathbb {R} ^ 1 \times \mathbb {R} ^ {n} \times \mathbb
{R} ^ {n \times n} $, called  {\em the second hypodifferential}
and {\em second hyperdifferential} respectively, for which the
representation $$ f (x + \Delta) = f (x) + \max_ {[a, v, A] \in
\underline {d} ^ 2 f (x)} [a + (v, \Delta) + \frac {1} {2} (A
\Delta, \Delta)] + $$ \be + \min_ {[b, w, B] \in \overline {d} ^ 2
f (x)} [b + (w, \Delta)) + \frac {1} {2} (B \Delta, \Delta)] + o_x
(\Delta ^ 2), \label {quasidiff2} \ee is true where $ \Delta ^ 2 =
\| \Delta \| ^ 2 $, $ o_x (\Delta) \rar 0 $, $ o_x (\Delta ^ 2)
\rar 0 $ with $ \Delta \rar 0 $, $ \frac {o_x (\al \Delta)} {\al}
\rar 0 $ and $ \frac {o_x (\al ^ 2 \Delta ^ 2)} {\al ^ 2} \rar 0 $
with $ \al \rar +0 $. The pairs of sets $ D  f (x) = [\underline
{d} f (x), \overline {d} f (x)] $ and $ D ^ 2 f (x) = [\underline
{d} ^ 2 f (x), \overline {d} ^ 2 f (x)] $ according to the
Demyanov's terminology are called {\em the first} and {\em the
second codifferentials} of  $ f $ at $ x $.  We suppose that the
function $ o_x (\Delta^2)$ in (\ref{quasidiff2}) is uniformly
infinitesimal with respect to $ \Delta^2  $ for small $ \Delta.  $

The author's problem is to define and construct the continuous
first and second codifferential of  $ f (\cdot) $ at $ x $, using
subdifferentials of the first and second orders, introduced by the
author in \cite{pimfirstsecondsubdiff}.

We use the quasidifferentials and codifferentials for writing the
necessary and sufficient conditions of optimality. It is necessary
to give some rules for their construction  which are just as
important as ability to calculate differentials for differentiable
functions.

Since the second-order subdifferential consists of symmetric
matrices of the second mixed derivatives of a function and we will
determine the second codifferential, using just such matrices, we
will assume further that the matrices, included in the second
codifferential, are symmetric matrices. It is not strong
limitation, because any matrix of the second codifferential can be
"symmetrized", i.e. be made symmetric, namely: instead of the
matrix $ A $ to consider symmetrized $ 1/2 (A + A ^ T), $ where $
A ^ T $ is transposed matrix. In this case, the equality (\ref
{quasidiff2}) is preserved.

To solve the formulated problem, the author used constructions
from \cite{pimfirstsecondsubdiff}. Let us recall them.

Further, to determine the subdifferentials of the first and second
order we need to remind the definition of the set-valued mappings
(SVM) $ D (\cdot) $ \cite{pimfirstsecondsubdiff}.

Consider a SVM $ D (\cdot) $ satisfying the following conditions:
\begin {enumerate}
\item $ x_0 \in \bint (x + D (x)) $ for all $ x \in S, S \in
\mathbb {R} ^ n  $ is a neighborhood of the point $ x_0 $;

\item The diameter of the set $ D (x) $, which we denote by $ diam
\, D (x) = d (D (x)), $ tends to zero as $ x \rar x_0 $, and
satisfies the inequality $ d (D (x)) \leq k \| x - x_0 \| $ for
some constant $k= k (D) $;

\item For some sequence $ \{\eps_i (D) \}, \eps_i \rar + 0, $ as $
i \rar \infty $ SVM $ D (\cdot) $ is constant for $ x $ from the
sets $ \eps_ {2i + 1} <\parallel x - x_0 \parallel <\eps_ {2i} $;

\item The boundary of the set $ D (x) $ for all $ x \in S, \, x
\neq x_0, $ is given by twice continuously differentiable
functions of $ x $. \end {enumerate} We will consider SVM $ D
(\cdot) $ satisfying written above conditions, for arbitrary
sequences $ \{\eps_i \}, \, \eps_i \rar + 0, $ and constants $ k
(D) $. Denote by $ \Xi $ the specified family of SVMs.

It is easy to give some examples of the considered family of SVMs.
We can take $D(x) =  B^n_{r(x)}=\{ w \in \mathbb{R}^n | \| w - x
\| \leq r(x) \},$ where $r(x)> \| x-x_0 \|$ and $r(x)  \rar 0$
when $x \rar x_0.$ The sets $D(x)$ can be the ellipsoids with the
similar qualities.

For an arbitrary Lipschitz function $ f (\cdot): \mathbb {R} ^ n
\rar \mathbb {R} $ we define the functions $ \varphi_{D} (\cdot):
\mathbb {R} ^ n \rar \mathbb {R} $ and $ \psi_{D} (\cdot): \mathbb
{R} ^ n \rar \mathbb {R} $
$$
\varphi_{D} (x): = \frac {1} {\mu (D (x))} \int_ {D (x)} f (x + y)
d y,
$$
$$
\psi_{D} (x): = \frac {1} {\mu (D (x))} \int_ {D (x)} \varphi (x +
y) d y,
$$
where $ D (\cdot): \mathbb {R} ^ n \Rightarrow {\mathbb {R} ^ n} $
is a SVM from the family of the SVM, defined above and in
\cite{proudintegapp1} - \cite{pimfirstsecondsubdiff}.

The set of curves was introduced in \cite{lowapp2}.
\begin{defi} $
\eta (x_0) $ is a set of smooth curves \mbox {$ r (x_0, \al, g) =
x_0 + \al g + o_r (\al) $}, where $ g \in S ^ {n-1} _1 (0) = \{v
\in \mathbb {R} ^ n: \| v \| = 1 \} $ and function \mbox {$ o_r
(\cdot): [0, \al_0] \rar \mathbb {R} ^ n, $} $ \al_0> 0, $
satisfies the following conditions:

\nid 1) $ o_r (\al) / \al $ tends to zero as $ \al \downarrow 0 $
uniformly in $ r (\cdot) $ ;

\nid 2) there is a continuous derivative $ o_r '(\cdot) $, and its
norm is bounded for all $ r $ in the following sense: there exists
$ c \; <\infty $ such that
$$
\sup _ {\tau \in (0, \al_0)} \parallel o_r '(\tau) \parallel \leq
c;
$$
3) the derivative $ \grad f (r (\cdot)) $ exists almost everywhere
(a.e.) along the curve $ r (x_0, \cdot, g) $. \label {defi1321}
\end{defi} We introduce the sets
$$
E f (x_0) = \{v \in \mathbb {R} ^ n: \exs \{{\al_ {k}} \}, \al_
{k} \downarrow 0, (\exs \, g \in S ^ {n-1} _1 (0)),
$$
$$ (\exs r (x_0, \cdot, g) \, \in \, \eta (x_0)),
v = \lim _ {\al_k \downarrow 0} \al_k ^ {- 1} \; \int ^ {\al_k} _0
\, \grad f (r (x_0, \tau, g)) d \tau \; \}
$$
and \be Df (x_0) = \bco \, \, Ef (x_0). \label {quasidiff2a} \ee

For the SVM $ D (\cdot) $ and the function $ f (\cdot) $ we
introduce the set
$$
\p \varphi_D (x_0) = \bco \{v \in \mathbb{R} ^ n \vl v = \lim_
{x_i \rar x_0} \varphi_ {D} '(x_i) \},
$$
where the points $ x_i $ are taken from the regions of constancy
of the SVM $ D (\cdot) $. The set $ \p \varphi_D (x_0) $ is convex
compact in $ \mathbb{R}^ n $ \cite{lowapp2a}.

Define the SVM $ \Phi f (\cdot): \mathbb {R} ^ n \rar 2 ^ {\mathbb
{R} ^ n} $ with images
$$
\Phi f (x_0) = \bco \, \bigcup_ {D (\cdot)} \, \p \varphi_D (x_0),
$$
where the union is taken over all the SVM $ D (\cdot) \in \Xi $. $
\Phi f (x_0) $ is called \cite{pimfirstsecondsubdiff} {\em the
subdifferential of the first order} of  $ f (\cdot) $ at the point
$ x_0 $. The theorem \cite{pimfirstsecondsubdiff} was proved,
establishing connection of the sets $ \Phi f (x_0) $ and $ Df
(x_0) $: for any Lipschitz function $ f (\cdot) $  equality
$$
\Phi f (x_0) = Df (x_0)
$$
is correct.

Also  {\em the second order subdifferential} was introduced in
\cite{pimfirstsecondsubdiff}.

We introduce the set of matrices for  $ f (\cdot) $
$$
\p ^ 2 \psi_D (x_0) = \bco \{A \in \mathbb {R} ^ {n \times n} \vl
A = \lim_ {x_i \rar x_0} \psi_ {D} '' (x_i) \},
$$
where the points $ x_i $ belong to the regions of constancy of the
SVM $ D (\cdot) \in \Xi $.

We define the SVM $ \Psi^2 f (\cdot): \mathbb {R}^ n \rar 2 ^
{\mathbb {R} ^ {n \times n}} $ with images
$$
\Psi ^ 2 f (x_0) = \bco \, \bigcup_ {D (\cdot)} \, \p ^ 2 \psi_D
(x_0),
$$
where the union is taken over all  SVM $ D (\cdot) \in \Xi $. We
call the set $ \Psi^2 f (x_0) $ by {\em the subdifferential of the
second order} of the function $ f (\cdot) $ at  $ x_0 $
\cite{pimfirstsecondsubdiff}.

We give an example of calculation of the subdifferentials.

It is known \cite{kutatelrubinov} that any finite convex  positive
homogeneous (p.h.) of the first order function $ h (\cdot):
\mathbb {R} ^ n \rar \mathbb {R} $ corresponds to a convex compact
set $\partial h (0) $ in $ {\mathbb {R} ^ n}, $ which is called
the subdifferential of the function $ h (\cdot) $ at $ 0 $. Also
the equality \be h (q) = \max_{v \in \partial h (0)} (v, q) \, \,
\, \, \forall q \in \mathbb{R}^ n, \label{quasidiff0} \ee is
correct, where $ (v, g)  $ is scalar product vectors $ v $ and $ g
$. Any convex compact set $ \partial h (0) $ corresponds to a
convex p.h. function $ h (\cdot) $ for which (\ref {quasidiff0})
is correct. This correspondence is called the Minkovsky duality.

We can find the set $ \partial h (0) $ corresponding to the
function $h (\cdot) $ in the following way. It is known that any
convex function is almost everywhere differentiable in $ \mathbb
{R} ^ n $. We denote the set, where the function $ h (\cdot) $ is
differentiable, by $ N_1 (h). $ Then
$$
\partial h (0) = \{v \in \mathbb {R} ^ n \vl \exists \{x_i
\}, x_i \rar 0, x_i \in N_1 (h), v = \lim_ {x_i \rar 0} h '(x_i)
\},
$$
which is called the Clarke subdifferential of the function $ h
(\cdot) $ at zero \cite{clark1}, \cite{clark}. We can calculate
$\partial h (0) $ in another way.

It was proved in \cite{lowapp2a} that for a convex  positively
homogeneous (p.h.) function $ h (\cdot): \mathbb {R} ^ n \rar
\mathbb {R} $
$$
\p h (0) = D h (0),
$$
as soon as $ \p h (0) $ coincides with the Clarke subdifferential
$ \p_ {CL} h (0). $

It was also proved in \cite{pimfirstsecondsubdiff}, that for a
convex p.h. of the first order function   $ h (\cdot);
\mathbb{R}^n \rar \mathbb {R} $
$$
\Phi h (0) = D h (0).
$$
So if
$$
h (q) = \max_{v \in B ^ n_1 (0)} (v, q),
$$
where $ B ^ n_1 (0) = \{v \in \mathbb {R} ^ n \vl \| v \| = 1 \} $
is the unit ball with center at zero in $ n- $ dimensional space,
then $ \Phi h (0) = D h (0) = \p h (0) = B ^ n_1 (0). $

\section {\bf Construction of a subdifferential
for a hypodifferentiable function}

Let $ f(\cdot): \mathbb {R} ^ n \rar \mathbb {R} $ be a Lipschitz
hypodifferentiable function at $ x $, i.e. the representation \be
f (x + \Delta) = f (x) + \max_ {[a, v] \in \underline {d} f (x)}
[a + (v, \Delta)]  + o_x (\Delta) \label{quasidiff3} \ee is true.
It is easy to get from  here
$$
f (x + \Delta) = f (x) + \max_ {v \in \underline \p f (x)} (v,
\Delta) + o_x (\Delta),
$$
where $ \underline \p f (x) = \{v \vl [\overline {a}, v] \in
\underline {d} f (x) \}, \overline {a} = \max \{a \vl [a, v] \in
\underline {d} f (x) \}. $ The set $ \underline \p f (x) $ is
called the subdifferential of the function $ f (\cdot) $ at a
point $ x $. Prove that $ \underline \p f (x) = D f (x). $

For any $ \Delta = \al g, \, \, g \in S ^ {n-1} _1 (0), $ we take
the curve $ r (x, \cdot, g) \in \eta (x), $ along which we
calculate the averaging integral of gradients
$$
\al ^ {- 1} \int ^ {\al} _0 \, \grad f (r (x, \tau, g)) d \tau.
$$
Take an arbitrary sequence $ \{\alpha_k \} $, $ \al_k \rar + 0, $
for which there is a limit
$$
v = \lim _ {\al_k \rar +0} \al_k ^ {- 1} \int ^ {\al_k} _0 \,
\grad f (r (x, \tau, g)) d \tau.
$$
By definition $ v \in D f (x). $ It is clear that
$$
(v, g) = \frac {\p f (x)} {\p g} = \max_ {u \in \underline \p f
(x)} (u, g).
$$
Therefore, $ D f (x) \supset \underline {\p} f (x). $

The strict inclusion is impossible, otherwise the vectors $ v \in
Df (x), v \notin \underline {\p} f (x), $ and $ g (v) \in S_1 ^ {n
- 1} (0): $
$$
g (v) = \frac {v - v_1} {\| v - v_1 \|},
$$
where $ v_1 = arg \min_ {w \in \p f (x)} \| v - w \| $, would
exist for which
$$
(v, g (v))> \frac {\p f (x)} {\p g (v)},
$$
which can not be true. This implies the equality $ Df (x) =
\underline {\p} f (x). $

So the following theorem is proven.
\begin{thm} If $ f (\cdot)  $ is a Lipschitz
hypodifferentiable function at a point $ x $, then $ D f (x) =
\underline \p f (x) $. \label{homogthm3}
\end{thm}


\section{Construction of the subdifferentials and
superdifferentials for codifferentiable functions}


Let $ f (\cdot): \mathbb {R} ^ n \rar \mathbb {R}  $ be a
Lipschitz codifferentiable function, i.e. the equality (\ref
{quasidiff1}) holds for it. It is not difficult to get from  (\ref
{quasidiff1}) \be f (x + \Delta) = f (x) + \max_ {v \in \underline
{\p} f (x)} (v, \Delta) + \min_ {w \in \overline {\p} f (x) } (w,
\Delta) + o_x (\Delta), \label{quasidiff6} \ee where
$$
\underline {\p} f (x) = \{v \vl [\bar {a}, v] \in \underline {d} f
(x) \}, \, \bar {a} = \max \{a \vl [a, v] \in \underline {d} f (x)
\}
$$
$$
\overline {\p} f (x) = \{w \vl [\bar {b}, w] \in \overline {d} f
(x) \}, \, \bar {b} = \min \{b \vl [b, w] \in \overline {d} f (x)
\}.
$$
The sets $ \underline {\p} f (x) $ and $ \overline {\p} f (x) $
are called {\em the subdifferential}  and {\em the
superdifferential} of the function $ f (\cdot) $ at the point $ x
$ correspondingly, and the function $ f (\cdot) $ is called {\em
quasidifferentiable} at $ x $.

Firstly, we consider the case when the function $ f (\cdot) $ can
be represented in the form
$$
f (x + \Delta) = f (x) + \max_ {v \in \underline {\p} f (x)} (v,
\Delta) + o_x (\Delta).
$$
Since the sets consisting from hypodifferentiable and
subdifferentiable functions coincide with each other
\cite{demrub}, it follows from Theorem \ref {homogthm3}

\begin{cor} If $ f (\cdot) $ is a Lipschitz subdifferentiable
function at $ x $, then $ Df (x) = \underline {\p} f (x). $ \label
{homogthm1} \end{cor} Introduce the difference of convex compact
sets $ A $ and $ B $ \cite{demrub}. \be A \rightharpoonup B =
\overline {co} \{\grad p_A (q) - \grad p_B (q) \vl q \in S ^ {n-1}
_A \cap S ^ {n-1} _B \}, \label {quasidiff7} \ee where $ p_A
(\cdot), p_B (\cdot) $ are the support functions to the sets $ A $
and $ B $ respectively:
$$
p_A (q) = \max_ {v \in A} (v, q), \, p_B (q) = \max_ {w \in B} (w,
q).
$$
$ S ^ {n-1} _A, \, S ^ {n-1} _B $ are the sets of the unit support
vectors to the sets $ A, \, B $, where the functions $ p_A
(\cdot), \, p_B (\cdot) $ are differentiable. The difference (\ref
{quasidiff7}) is called the {\em Demianov difference}
\cite{demrub}.

Couples of convex compact sets $ [A, B] $ and $ [C, D] $ are
called {\em equivalent} if $ A \rightharpoonup B = C
\rightharpoonup D $.

Let the equality (\ref{quasidiff6}) be true. Rewrite (\ref
{quasidiff1}) in the form (\ref{quasidiff6}). Take $ \Delta = \al
g, \, g \in S ^ {n-1} _1 (0)$  for which there is a sequence $
\{\al_k \}, \al_k \rar + 0, $ that the following equalities are
true for the vector $ g $ and a vector $ v \in Df (x):$
$$
v = \lim _ {\al_k \rar +0} \al_k ^ {-1} \; \int ^ {\al_k} _0 \,
\grad f (r (x, \tau, g)) d \tau,
$$
$$
(v, g) = (v_1, g) - (v_2, g) = \frac {\p f (x)} {\p g},
$$
where
$$
(v_1, g) = \max_ {w \in \underline{\p} f (x)} (w, g), \, \, (v_2,
g) = \max_ {w \in - \overline {\p} f (x)} (w, g).
$$
It follows that
$$
\underline{\p} f (x) \rightharpoonup (- \overline{\p} f (x) )
\subset Df (x).
$$
Strict inclusion can not be, otherwise the vectors $ v \in Df (x)
$ and $ g (v) \in S ^ {n-1} _1 (0) $ would exist, for which
$$
\frac {\p f (x)} {\p g (v)} = (v, g (v))> (v_1, g (v)) - (v_2, g
(v)) = \frac {\p f (x)} {\p g (v)},
$$
where
$$
(v_1, g (v)) = \max_ {w \in \underline {\p} f (x)} (w, g (v)), \,
(v_2, g (v)) = \max_ {w \in - \overline {\p} f (x)} (w, g (v),
$$
what can not be correct. So the following theorem is proven.
\begin{thm} If $ f (\cdot)$ is a
Lipschitz quasidifferentiable function at a point $ x, $ then
$$
Df (x) = \underline {\p} f (x) \rightharpoonup (- \overline {\p} f
(x)).
$$
The sets $ \underline {\p} f (x) $ and $ \overline {\p} f (x) $
are determined by the set $ Df (x) $ up to equivalence.
\label{homogthm2} \end{thm}

We will start to calculate the second codifferentials of  p.h.
functions of the second order.

\section {\bf Positively homogeneous functions of the second
order}

Now let $ h (\cdot): \mathbb {R} ^ n \rar \mathbb {R} $ be p.h. of
the second-order function that is twice hypo-differentiable at
zero, i.e.
$$
h (\lambda q) = \lambda ^ 2 h (q) \, \, \, \, \forall \lambda> 0,
\forall q \in \mathbb {R} ^ n
$$
and the equality \be h (q) = \max_{[a, v, A] \in \underline{d}^2
h(0)} [a + (v, q) + \frac {1} {2} (A q, q)] + o (q ^ 2), \label
{quasidiff4} \ee is true  where $ \frac {o (\al ^ 2 q ^ 2)} {\al ^
2} \rar 0 $ as $ \al \rar +0 $. The set $ \underline {d} ^ 2 h (x)
$ is called {\em the second hypodifferential} according to the
terminology by V.F. Demyanov.

We now rewrite (\ref{quasidiff4}) in the form
$$
h (q) = \max_ {[v, A] \in \underline {\p} ^ 2 h (0)} [(v, q) +
\frac {1} {2} (A q, q)] + o (q ^ 2),
$$
where $ \underline {\p} ^ 2 h (0) = \{[v, A] \vl \exists \,
[\overline {a}, v, A] \in \underline {d} ^ 2 h (0) \}, \overline
{a} = \max \{a \vl \exists \, [a, v, A] \in \underline {\p} ^ 2 h
(0) \}. $ The set $ \underline {\p} ^ 2 h (0) $ is called {\em the
second subdifferential} of the function $ h (\cdot) $ at zero.

By the assumption $ h (\cdot) $ is a p. h. function of the  second
order. Therefore terms in the expansion (\ref {quasidiff4}),
containing the linear functions with respect to $ q $, will not
be, i.e. \be h (q) = \max_ {A \in \cal A} \, \frac {1} {2} (Aq,
q), \label {quasidiff5} \ee where $ {\cal A} = \{A [n \times n]
\vl \exists \, [v, A] \in \underline {\p} ^ 2 h (0) \} $. The
following theorem is proven.

\begin{thm} Let $ h (\cdot): \mathbb {R} ^ n \rar \mathbb {R} $
be a twice hypo-differentiable p.h. function of the second order,
then there exists a convex, compact set of matrices $ \cal {A} $
such that the equality (\ref{quasidiff5}) is true. \end{thm}

We will prove that the function $ h (\cdot) $ is a.e. twice
differentiable in $ \mathbb{R}^ n $. Note that in the case of
convexity of the function $ h (\cdot) $ the last statement follows
from the well known Alexandrov's theorems \cite{alexandrovAD},
which states that any finite convex  function is a.e. twice
differentiable in $\mathbb{R}^n$.

Let us take an arbitrary vector $ \bar {q} \in S_1 ^ {n-1} (0). $
Denote by
$$
V (\bar {q}) = \{\bar {A} \in {\cal A} \vl (\bar {A} \bar {q},
\bar {q}) = \max_ {A \in {\cal A}} (A \bar {q}, \bar {q}) \}.
$$
If $ V (\bar {q}) $ consists of a single matrix $ \bar {A} $, then
the function $ h (\cdot) $ is twice differentiable at the point  $
\bar {q} $ and $ h''(\bar {q}) = \nabla ^ 2 h (\bar {q}) = \bar
{A}.$

The last statement can be proved in the same way how it was done
for the differentiability of a  function, represented in the form
(\ref{quasidiff0}), at points $ q $, where the maximum by $ v \in
\partial h (0) $ is reached at a single point. In our case the
function (\ref{quasidiff5}) can be represented as the maximum of
the scalar product of vectors whose coordinates are expressed by
the elements of the matrix $ A $ and the coordinates of the vector
$ q $.

It follows from the convexity of the set $ \cal {A} $ and the
above that the function $ h (\cdot) $ is a.e. twice differentiable
in $ \mathbb {R} ^ n $. Denote the set, where the function $ h $
is twice differentiable in $ \mathbb {R} ^ n $, by $ N_2 (h) $.

Let us prove that
$$
{\cal A} = \bco \{A [n \times n] \vl A = \lim_ {q_i \rar 0} h ''
(q_i), \, \, q_i \in N_2 (h) \}.
$$

The set $ \cal A $ can be approximated in the Hausdorff metrics
with any precision by the set $ {\cal A} _m $ with a smooth
boundary so, what a set
$$
V_m (\bar {q}) = \{\bar {A} \in {\cal A} _m \vl (\bar {A} \bar
{q}, \bar {q}) = \max_ {A \in {\cal A} _m} (A \bar {q}, \bar {q})
\}
$$
will consist of one matrix $ \bar A $ for any vector $ \bar q. $
The last statement means that the function
$$
h_m (q) = \max_ {A \in {\cal A} _m} \, \frac {1} {2} (Aq, q),
$$
is twice differentiable on $ B ^ n_1 (0) \backslash \{0 \}, $
where $ B ^ n_1 (0) = \{v \in \mathbb {R} ^ n \vl \| v \| \leq 1
\}. $ Such approximation guarantees us an uniform approximation of
matrices $ \nabla ^ 2 h (q) $   with the help of matrices $ \nabla
^ 2 h_m (q) $ on $ B ^ n_1 (0) \backslash \{0 \} $. But
$$
{\cal A} _m = \bco \{A_m [n \times n] \vl A_m = \lim_ {q_i \rar 0}
h_m '' (q_i), \, \, q_i \in B ^ n_1 (0) \backslash \{0 \} \}.
$$
It follows from the above that
$$
{\cal A} = \bco \{A [n \times n] \vl A = \lim_ {q_i \rar 0}
h"(q_i), \, \, q_i \in N_2 (h) \},
$$
as soon as
$$
\lim_ {m \rar \infty} \rho_H ({\cal A} _m, {\cal A}) = 0,
$$
where $ \rho_H  $ is the Hausdorff metrics.

We get the following result.
\begin{thm} If  the equality (\ref
{quasidiff5}) is true for a function $ h (\cdot): \mathbb {R} ^ n
\rar \mathbb {R} $, then
$$ {\cal
A} = \bco \{A [n \times n] \vl A = \lim_ {q_i \rar 0} h '' (q_i),
\, \, q_i \in N_2 (h) \}.
$$
\end{thm}

We will show now that $ {\cal A} = \Psi ^ 2 h (0), $ i.e. $ \cal A
$ coincides with the second-order subdifferential of the function
$ h (\cdot) $ at the point zero.

At all points $ q $, where the function $ h (\cdot) $ is twice
differentiable, $ \Psi ^ 2 h (q) = \{h '' (q) \}. $ But if the
function $ h (\cdot) $ is twice differentiable at $ \bar {q}, $,
then it is twice differentiable at all points of the ray $
\{\lambda \, \bar {q} \vl \lambda> 0 \} $, moreover, $ h ''
(\lambda \bar q) = \bar A. $ With a suitable choice of the
set-valued mapping (SVM) $ D (\cdot) $ from the definition of the
second order subdifferential \cite {pimfirstsecondsubdiff} we
obtain a matrix $ A $ arbitrarily close to the matrix $ \bar A =
h''( \bar q). $ Therefore, $ {\cal A} \subset \Psi ^ 2 h (0). $

On the other hand, $ \Psi ^ 2 h (0) $ is, by definition, the
convex envelope of all limit matrices $ h '' (q_i) $ as $ q_i \rar
0. $ Hence, $ \Psi ^ 2 h (0) \subset {\cal A}. $ The two
inclusions imply the equality $ \Psi ^ 2 h (0) = {\cal A}. $

So the following theorem is proven.

\begin{thm} If  the equality (\ref {quasidiff5}) is true
for a function $ h (\cdot): \mathbb {R} ^ n \rar \mathbb {R} $,
then $ \Psi ^ 2 h (0) = {\cal A}, $ i.e. the set  $ {\cal A} $
coincides with the subdifferential of the second order of the
function $ h (\cdot) $ at $ 0 $. \end{thm}

\section {\bf The second codifferentials}

Let $ f (\cdot): \mathbb {R} ^ n \rar \mathbb {R} $ be a Lipschitz
twice codifferentiable function at $ x $, i.e. the equality (\ref
{quasidiff2}) is true. Let us consider from the beginning the case
when $ f (\cdot) $ is twice hypo-differentiable, i.e. \be f (x +
\Delta) = f (x) + \max_ {[a, v, A] \in \underline {d} ^ 2 f (x)}
[a + (v, \Delta) + \frac {1} {2} (A \Delta, \Delta)] + o_x (\Delta
^ 2). \label {quasidiff8} \ee We can get from (\ref {quasidiff8})
the following equality \be f (x + \Delta) = f (x) + \max_ {[v, A]
\in \underline {\p} ^ 2 f (x)} [(v, \Delta) + \frac {1} {2} (A
\Delta, \Delta)] + o_x (\Delta ^ 2), \label {quasidiff9} \ee where
$$
\underline {\p} ^ 2 f (x) = \{[v, A] \vl \exists \, [\overline
{a}, v, A] \in \underline {d} ^ 2 f (x) \}, \, \overline {a} =
\max \{a \vl \exists \, [a, v, A] \in \underline {d} ^ 2 f (x) \}.
$$
The set $ \underline {\p} ^ 2 f (x) $ is called {\em the  second
subdifferential} of the function $ f (\cdot) $ at the point $ x $,
which is a convex compact set, that follows from the properties of
the set $ \underline {d} ^ 2 f (x). $

Let us define the set
$$
\underline {\p} f (x) = \{v \vl \exists \, [v, A] \in \underline
{\p} ^ 2 f (x) \},
$$
which is called the {\em subdifferential} of $ f (\cdot)$ at $x$.
It is a convex compact set, that follows from the convexity and
compactness of the sets $ \underline {d} ^ 2 f (x) $ and $
\underline {\p} ^ 2 f (x) $. We get the decomposition from (\ref
{quasidiff8}) and (\ref {quasidiff9}) for an arbitrary $ \triangle
= \al g, \al> 0, $
$$
f (x + \al g) = f (x) + \al \max_ {v \in \underline {\p} f (x)}
(v, g) + o (\al).
$$
This implies the equality for the directional derivative
$$
\frac {\p f (x)} {\p g} = \max_ {v \in \underline {\p} f (x)} (v,
g).
$$
We introduce a function $ \tilde {f} (\cdot): \mathbb {R} ^ n \rar
\mathbb {R} $ \be \tilde {f} (x + \triangle) = f (x + \triangle) -
\max_{v \in \underline {\p} f (x)} (v, \triangle) - f (x) = \frac
{1} {2} \max_ {A \in {\cal A}} (A \triangle, \triangle) + o_x
(\triangle ^ 2), \label {quasidiff10} \ee where $ {\cal A} = \{A
\vl \exists [v, A] \in \underline {\p} ^ 2 f (x), \, v \in
\underline {\p} f (x) \}. $ The set $ \cal A $ is convex and
compact, which follows from the convexity and compactness of the
sets $ \underline {d} ^ 2 f (x) $ and $ \underline {\p} ^ 2 f (x)
$.

The validity of the representation (\ref {quasidiff10}) follows
from (\ref {quasidiff9}). If (\ref {quasidiff10}) were wrong,
there would not be a pair $ [v, A] \in \underline {\p} ^ 2 f (x)
$, for which the equality (\ref {quasidiff9}) was true for all
small $ \Delta $.

We can prove in the same way, as it was done for the p.h.
functions of the second order, that
$$
\Psi ^ 2 \tilde {f} (x) = {\cal A},
$$
because of the additive $ o_x (\triangle ^ 2) $ is not important
for calculation of the second order subdifferential. This fact can
be proved in same way as it was done in
\cite{pimfirstsecondsubdiff} for a twice differentiable function
at the point $ x $.  It was proved in \cite
{pimfirstsecondsubdiff} for this case that the subdifferential of
the second order coincides with the second derivative of this
function at $x$.

Thus, the algorithm for calculation of the second subdifferential
$ \underline {\p} ^ 2 f (x) $ of the function $ f (\cdot) $ at $ x
$ in the case of its representations in the form (\ref
{quasidiff8}), i.e. when $ f (\cdot) $ is hypo-differentiable, is
following:

1) we represent the function $ f (\cdot) $ in the form (\ref
{quasidiff9});

2) we find the subdifferential $ \underline {\p} f (x) $ according
to the theorem \ref {homogthm1};

3) we find the function $ \tilde {f} (\cdot) $ by the formula
(\ref {quasidiff10});

4) we find the second-order subdifferential of the function $
\tilde {f} (\cdot) $ at $ x $, i.e. find the set $ \Psi ^ 2 \tilde
{f} (x), $ which coincides with $ \cal A $.

The second hypo-differential for a hypo-differentiable function $
f (\cdot) $ will be equal to
$$
\underline {d} ^ 2 f (x) = \bco \{[a, v, A] \vl v \in \underline
{\p} f (x)= Df(x), \, A \in \Psi ^ 2 \tilde {f} (x) = {\cal A} \},
$$
where $ a \in [-a_0, 0], a_0> 0, $ and $ [0, v, A] \in \underline
{d} ^ 2 f (x) $ only for $ v, A $  belonging  to the boundaries of
the sets $ \underline {\p} f (x) $ and $ \Psi ^ 2 \tilde {f} (x) =
{\cal A} $ respectively, or, in other words, to the boundary of
the second subdifferential $ \underline {\p} ^ 2 f (x) $.

Let now the function $ f (\cdot) $ be twice codifferentiable at
the point $ x $, i.e. the equality (\ref {quasidiff2}) holds.
Rewrite (\ref {quasidiff2}) in the form
$$
f (x + \Delta) = f (x) + \max_ {[v, A] \in \underline {\p} ^ 2 f
(x)} [(v, \Delta) + \frac {1} {2} (A \Delta, \Delta)] +
$$
$$
+ \min_ {[w, B] \in \overline {\p} ^ 2 f (x)} [(w, \Delta)) +
\frac {1} {2} (B \Delta, \Delta)] + o_x (\Delta ^ 2),
$$
where $ \underline {\p} ^ 2 f (x), \overline {\p} ^ 2 f (x)$ is
the second subdifferential and superdifferential respectively, $
[\underline {\p} ^ 2 f (x), \overline {\p} ^ 2 f (x)] $ is the
second quasidifferential of the function $ f $ at the point $ x $.

We can take any vectors $ v, w $ from the sets $ \underline {\p} f
(x), \overline {\p} f (x) $ respectively, called the
subdifferential and superdifferential, for which (\ref
{quasidiff6}) is true. The sets $ \underline {\p} f (x), \overline
{\p} f (x) $ are defined with accuracy to equivalence \be Df (x) =
\underline {\p} f (x) \rightharpoonup (- \overline {\p} f (x)).
\label {quasidiff11} \ee Let the the sets $ \underline {\p} f (x),
\overline {\p} f (x) $ be already defined. We introduce the
function
$$
\tilde {f} (x + \Delta) = f (x + \Delta) - f (x) - \max_ {v \in
\underline {\p} f (x)} (v, \Delta) - \min_ {w \in \overline {\p} f
(x) } (w, \Delta) =
$$
$$
= \frac {1} {2} \max_ {A \in {\cal A}} (A \triangle, \triangle) +
\frac {1} {2} \min_ {B \in {\cal B}} (B \triangle, \triangle) +
o_x ( \triangle ^ 2).
$$
According to the previous arguments we have

\begin{thm}
The matrices $ A, B $ can be taken from the sets $ {\cal A} $ and
$ {\cal B} $, for which \be \Psi ^ 2 \tilde {f} (x) = {\cal A}
\rightharpoonup {(- \cal B)}, \label {quasidiff12} \ee
where
$$
{\cal A} = \{A \vl \exists [v, A] \in \underline {\p} ^ 2 f (x),
\, v \in \underline {\p} f (x) \}, \, \, {\cal B} = \{B \vl
\exists [w, B] \in \overline {\p} ^ 2 f (x), \, w \in \overline
{\p} f (x) \}.
$$
\end{thm}

Thus, the sets $ {\cal A} $ and $ {\cal B} $ also determined
ambiguously, namely: accurately to the equivalence so that the
equality (\ref {quasidiff12}) would hold.

In \cite {demrub} on page 189 the authors write: "Of course, the
answer for the question of how constructively to find the set $
\underline {d} f (x) $ for a convex function is significant for
practical use. This question must be solved for specific classes
of convex functions."  Let us answer this question.

If $ f (\cdot)  $ is convex finite in $ \mathbb {R} ^ n $, then it
is hypodifferentiable and the representation \be f (x + \Delta) =
f (x) + \max_ {[a, v] \in \underline {d} f (x)} [a + (v, \Delta)]
\label {quasidiff13} \ee is true. It is not difficult to move from
(\ref {quasidiff13}) to the form
$$
f (x + \Delta) = f (x) + \max_ {v \in \underline \p f (x)} (v,
\Delta) + o (\Delta),
$$
where $ \underline \p f (x) = \{v \vl [\overline {a}, v] \in
\underline {d} f (x) \}, \overline {a} = \max \{a \vl [a, v] \in
\underline {d} f (x) \}. $

It was shown (Theorem \ref {homogthm1}), that $ Df (x) =
\underline {\p} f (x) $. We can take as a segment of values of the
parameter $ a $ any segment $ [- a_0, 0] $, $ a_0> 0, $ for which
$ [0, \bar {v}] \in \underline {d} f (x) $ only for $ \bar {v} $
belonging to the boundary of the set $ \underline \p f (x) $.

The above is true for an arbitrary hypodifferentiable functions.



The construction of the second codifferential $ D ^ 2 f (x) =
[\underline {d} ^ 2 f (x), \overline {d} ^ 2 f (x)] $ at  $ x $
can be done in the same way as the construction of the first
codifferential. The second hypodifferential and hyperdifferential
for a codifferentiable function $ f (\cdot) $  are equal to
$$
\underline {d} ^ 2 f (x) = \bco \{[a, v, A] \vl v \in \underline
{\p} f (x), \, A \in {\cal A} \},
$$$$
\overline {d} ^ 2 f (x) = \bco \{[b, w, B] \vl w \in \overline
{\p} f (x), \, B \in {\cal B} \},
$$
where $ a \in [-a_0, 0], a_0> 0, $ and $ [0, v, A] \in \underline
{d} ^ 2 f (x) $ only for $ v, A $, belonging to the boundary  of
the second subdifferential $ \underline {\p} ^ 2 f (x) $, $ b \in
[0, b_0], b_0
> 0, $ and $ [0, w, B] \in \underline {d} ^ 2 f (x) $ only for $ w, B $,
belonging to the boundary of the second superdifferential $
\overline {\p} ^ 2 f (x). $

\section {Construction of a continuous first
codifferential}

Although we have shown how to build the first and second
codifferentials for Lipschitz functions, but the constructed
codifferentials are not necessarily  continuous set-valued
mappings (SVM), as functions of $ x $. The main goal of the
introduction of the codifferentiable functions was precisely, that
the codifferentials were continuous SVM in contrast to the
predecessors of them, namely, the quasidifferentials. Our further
problem  is construction the first and second continuous
codifferentials for the Lipschitz functions. At first we will
limit ourselves to building continuous first codifferentials.
Further  we  will consider how to the construct the second
continuous codifferentials.

It follows from the above that it is sufficiently to confine
consideration for the hypo-differentiable functions, i.e. for
functions for which  the decomposition
$$
f (x + \Delta) = f (x) + \max_ {[a, v] \in \underline {d} f (x)}
[a + (v, \Delta)] + o_x (\Delta)
$$
is true. Construction of the  first continuous codifferential for
the Lipschitz function $ f (\cdot) $ will be based on the
following theorem.

\begin{thm} Any subdifferentiable Lipschitz function
$ f (\cdot): \mathbb {R} ^ n \rar \mathbb {R} $, whose
subdifferential $ \p f(\cdot) $ is upper semi-continuous
(\cite{demrub}, p. 399) as the SVM at $ x $, is continuously
codifferentiable at $ x $. \label{homogthm4} \end{thm} {\bf
Proof.} For any $ g \in S_1 ^ {n-1} (0) = \{w \mid \parallel w
\parallel = 1 \} $ decomposition takes place
$$
f (x + \al g) = f (x) + \al \max_ {v \in \p f (x)} (v, g] + o
(\al, g) = f (x) + \al f '(x, g) + o (\al, g),
$$
where by definition
$$
f '(x, g) = \frac {\p f (x)} {\p g}, \, \, \, \, \lim _ {\al \rar
+0} \frac {o (\al, g)} {\al} = 0.
$$
Note that from the boundedness of the direction derivative of the
function $ f (\cdot) $ in a neighborhood of an arbitrary point and
upper semicontinuity (UPSC) of  $ \p f (\cdot) $ it follows the
Lipschitz quality of the function $ f (\cdot) $  with a Lipschitz
constant depending on properties of the mapping $ \p f (\cdot) $.

It was proved (Theorem \ref {homogthm3}) that
$$
D f (x) = \p f (x),
$$
where the set $ D f (x) $ was defined earlier (see
(\ref{quasidiff2a})).

Let us define for some $ \al_0> 0 $ the set
$$
V (\al, x) = \overline {\bco} \{v \in \mathbb {R} ^ n \mid \exists
g \in S_1 ^ {n-1} (0), \exists r (x, \cdot, g) \in \eta (x), v =
\al ^ {- 1} \; \int ^ {\al} _0 \, \grad f (r (x, \tau, g)) d \tau
\},
$$
$$
V (x) = \overline {\bco} \bigcup _ {\al \in (0, \al_0]} V (\al,
x).
$$
The set $ V (x)$ is bounded, which follows from the Lipschitz
quality of the function $ f (\cdot) $.

From the equality
$$
V (0, x) = Df (x) = \p f (x),
$$
it follows that SVM $ V (\cdot): \mathbb {R} ^ n \rar 2 ^ {\mathbb
{R} ^ n} $  is a continuous extension of the subdifferential
mapping $ \p f (\cdot) $ (see \cite {lowapp2}).

Construct the hypodifferential of the function $ f (\cdot) $ at
the point $ x $
$$
\underline {d} f (x) = \bco \{[a, v] \in \mathbb {R} ^ {n + 1}
\mid \exists \beta \in (0, \al_0], \exists g \in S_1 ^ {n-1} (0),
\exists r (x, \cdot, g) \in \eta (x),
$$
$$
a (\beta, g) = - \beta \rho (v (\beta, g), \p f (x)), \, \, v
(\beta, g) = \beta ^ {- 1} \; \int ^ {\beta} _0 \, \grad f (r (x,
\tau, g)) d \tau \} \}.
$$
Here
$$
\rho (v (\beta, g),\p f (x)) =  \bmin_ {w \in \p f (x)} \| v
(\beta, g) - w \|
$$
is the deviation of $v (\beta, g)$ from $\p f (x)$.

Check the expansion for $ \al> 0 $ and $ g \in S_1 ^ {n-1} (0) $
\be f (x + \al g) = f (x) + \max _ {[a (v), v] \in \underline {d}
f (x)} [a + \al (v, g)] + o (\al, g). \label {quasidiff15} \ee
Note that \be f (x + \al g) -f (x) = \al (\al ^ {- 1} \int ^ {\al}
_0 \, \grad f (r (x, \tau, g)) d \tau, g) = \al (v (\al, g), g),
\label {quasidiff16} \ee where
$$
v (\al, g) = \al ^ {- 1} \int ^ {\al} _0 \, \grad f (r (x, \tau,
g)) d \tau.
$$
Denote by
$$
\overline {v} (g) = arg \max_ {w \in \p f (x)} (w, g).
$$
As soon as
$$
\lim_{\al \rar +0} \, \rho (v (\al, g), \p f (x)) = 0,
$$
then the equality (\ref {quasidiff16}) can be rewritten as
$$
f (x + \al g) -f (x) = \al (\overline {v} (g), g) + \al (v (\al,
g) - \overline {v} (g), g) = \al (\overline {v} (g), g) + o (\al,
g) =
$$
\be = \al (v (\al, g), g) - \al \rho (v (\al, g), \p f (x)) +
\overline {o} (\al, g), \label {quasidiff17} \ee where by the
definition
$$
\lim_{\al \rar +0} \frac {\overline {o} (\al, g)} {\al} = 0.
$$
We will show that for all $ \beta \in (0, \al_0] $
$$
f (x + \al g) \meq f (x) + a (\beta, g) + \al (v (\beta, g), g) +
o (\al, g).
$$
The following inequality
$$
f (x) + \al (v (\beta, g), g) - \beta \rho (v (\beta, g), \p f
(x)) \leq f (x) + \al (\overline {v} (g), g) +
$$
\be + \al \max_ {v (\beta, g) \in V (x)} (v (\beta, g) - \overline
{v} (g), g) - \beta \rho (v (\beta, g), \p f (x)). \label
{quasidiff18} \ee is correct. Since for any $ \beta \in (0, \al_0]
$
$$
\rho (v (\beta, g), \p f (x)) \meq (v (\beta, g) - \overline {v}
(g), g),
$$
then  the inequality
$$
f (x) + \al (v (\beta, g), g) - \beta \rho (v (\beta, g), \p f
(x)) \leq
$$
\be \leq f (x) + \al (\overline {v} (g), g) + o (\al, g) = f (x +
\al g) \label {quasidiff19} \ee follows from (\ref{quasidiff18})
for small $ \al $ and all $ \beta> \al. $

As it follows from (\ref{quasidiff17}), the inequality
(\ref{quasidiff19}) turns into equality for $ \beta = \al. $ For $
\beta <\al $
$$
\al \max_ {v (\beta, g) \in V (x)} (v (\beta, g) - \overline {v}
(g), g) - \beta \rho (v (\beta, g), \p f (x)) = \tilde {o} (\al,
g),
$$
as soon as
$$
\lim _ {\beta \rar 0} (v (\beta, g) - \overline {v} (g), g) = 0,
$$
$$
\lim _ {\beta \rar 0} \rho (v (\beta, g), \p f (x)) = 0,
$$
where
$$
\lim _ {\al \rar +0} \frac {\tilde {o} (\al, g)} {\al} = 0.
$$
Therefore we rewrite (\ref {quasidiff18}) as \be f (x) + \al (v
(\beta, g), g) - \beta \rho (v (\beta, g), \p f (x)) \leq f (x) +
\al (\overline {v} (g), g) + \tilde {o} (\al, g) = f (x + \al g)
\label{quasidiff20} \ee Note that in the decomposition
(\ref{quasidiff20}) the function $ \tilde{o}(\cdot, g)$ is
uniformly infinitely small in $ g \in S_1^{n-1}(0)$, since in the
decomposition
$$
f (x + \al g) = f (x) + \al f '(x, g) + o (\al, g)
$$
the function $ o (\al, g) $ is uniformly infinitely small in $ g
\in S_1 ^ {n-1} (0) $, which follows from the Lipschitz quality of
$ f (\cdot) $  and the upper semicontinuity of $ \p f (\cdot) $.

So, we considered all possible cases and proved the equality (\ref
{quasidiff15}). The theorem is proved. $ \Box $

We will formulate the theorem the proof of which follows from the
theorem \ref{homogthm4}.

\begin{thm} Any quasidifferentiable Lipschitz function $ f (\cdot):
\mathbb {R} ^ n \rar \mathbb {R} $, the SVMs $\underline {\p} f
(\cdot)$, $\overline {\p} f (\cdot)$ of which  are upper
semi-continuous at $ x $ and, consequently, the SVM $ D f (\cdot)
$ is upper semicontinuous, is continuously codifferentiable at $ x
$. \label{homogthm5}
\end{thm}

Let us give an example of a Lipschitz quasidifferentiable function
$ f (\cdot) $ for which $ Df (\cdot) $ is not upper semicontinuous
and the function $ f (\cdot) $ is not continuously
codifferentiable.

\begin {ex} \hspace {-2mm}. The graph of the function $ f: \mathbb
{R} \rar \mathbb {R} $ consists of the segments,  located between
the curves $ -x ^ 2, + x ^ 2 $, with the slopes $ \pm 1 $. The
function $ f (\cdot) $  is not representable as the difference of
two convex functions in a neighborhood of the point zero, since $
\vee ^ a _0 f '= \infty $ for an arbitrary $ a
> 0. $ It is easy to see that $ \p _ {Cl} f (0) = [-1, + 1], \,
D f (0) = \{0 \}.$ The SVM $ Df (\cdot) $ is not upper
semi-continuous at zero. The function $ f (\cdot) $ is
quasidifferentiable, but not continuously codifferentiable at
zero. \label{ex1322} \end{ex}

\section{Applications}

Let us consider the problem of finding a minimum of a Lipschitz
quasidifferentiable function $f(\cdot)$. The quasidifferentiable
functions are the generalization of convex functions and any ones
that are the difference of convex functions. The application of
these functions is wide. The methods of optimization of them are
described in huge list of publications, which can be found in
\cite{demvas1}, \cite{demrub}, \cite{demyanovdixon} -
\cite{vygodchikova}. The methods of optimization of a
quasidifferentiable function $f(\cdot)$ \cite{demvas1} are based
on calculations of the subdifferentials and superdifferentials and
the necessary  conditions for minimum at $x^* \in \mathbb{R}^n$
\cite{demvas1},\cite{demrub}
$$
- \overline {\p} f (x^*) \subset  \underline {\p} f (x^*).
$$
For a subdifferentiable function $f(\cdot)$ the above conditions
are written in the following form
$$
0 \in {\p} f (x^*).
$$
We can conclude from Theorem \ref{homogthm2} the following
\begin{thm}
The inclusion
$$
0 \in Df(x^*)
$$
is the necessary condition of optimality.
\end{thm}

If the inclusion \be - \overline {\p} f (x^*) \subset  \bint
\underline {\p} f (x^*) \label{quasidiff21} \ee is correct, then
the point $x^*$ is an minimal point i.e. (\ref{quasidiff21}) is an
sufficient condition of minimum at $ x^*$ \cite{demvas1}.

Basing on the result of Theorem \ref{homogthm2} and the written
above, we have
\begin{cor}
If
$$
0 \in \bint Df(x^*),
$$
then $x^*$ is an optimal point of a quasidifferentiable function
$f(\cdot)$.
\end{cor}

Without knowing the subdifferential and superdifferential we can
not apply the optimization methods. Therefore, the methods of
finding for the subdifferentials and superdifferentials are as
important for practical optimization of quasidifferentiable
Lipschitz functions as the rules of calculations of derivatives
for differentiable functions. The subdifferentials of the second
order are necessary for development of methods of the second order
for nonsmooth functions.

Let us provide the sufficient minimum conditions which are similar
to the sufficient minimum conditions for smooth functions
\cite{pimfirstsecondsubdiff}.

\begin{thm}
If the necessary condition for the minimum of $f(\cdot)$ at $x^*$
is true and there exists $\beta(g) >0$ for all directions $g \in
G$ for which $f'(x^*,g)=\frac{\p f(x^*)}{\p g}=0$ and the
inequality
$$
(Ag,g) \meq \beta(g) \| g  \|^2 \,\,\,\,\,\, \forall A \in
\p^2\psi_D (x^*), \,\,\, \forall D(\cdot) \in \Xi
$$
holds, then $x^*$ is the minimum of $f(\cdot)$.
\label{thmfirstsecsub6}
\end{thm}

Let us give an example proving importance of the giving
constructions.

\begin{ex}
Let be $f(\cdot): \mathbb{R} \rar \mathbb{R}$ with a graph lying
between two curves $y= \vl x \vl + x^2 $ and $y= \vl x \vl - x^2 $
and consisting from slopes $\pm 1$ with  the limit point  at zero.
Then
$$
Df(0)=[ -1; 1 ].
$$
We can conclude from $0 \in \bint Df(0) $  that the point zero is
an optimal point.
\end{ex}

\begin{ex}
Let be $f(\cdot): \mathbb{R} \rar \mathbb{R}$ with a graph lying
between two curves $y=x^2$ and $y=2x^2$ and consisting from slopes
$\pm 1$ with  the limit point  at zero. Then
$$
Df(0)=\{ 0 \}, \,\,\,\, \Psi^2 f(0)= [2 ; 4].
$$
We can conclude from here that zero is a minimum point.
\end{ex}

To optimize the function from the Example 1 we can use the
function $\varphi_D (\cdot) $ using the gradient methods with a
SVM $D(\cdot)$. To optimize the function from the Example 2 we can
use the function $\psi_D (\cdot) $ using the gradient methods with
a SVM $D(\cdot)$. During the optimization process we reduce the
diameter $d(D(x_k))$ of $D(x_k)$. We make consistent reduction of
$d(D(x_k))$ with respect to the length of step $ \| \Delta_k \| $.
The details of the optimization process were written in another
article.

\section{Conclusion}

In the article the author is studying  the twice codifferentiable
functions, introduced by Prof. V.Ph. Demyanov,  and some methods
for calculating their codifferentials.   At the beginning the
easier case is considered  when a function is twice
hypodifferentiable. There is proved that a twice
hypodifferentiable positively homogeneous function $ h (\cdot) $
of the second order is maximum of the quadratic forms with respect
to a certain set of matrices, which coincides with the convex hull
of the limit matrices calculated at points, where the original
function $ h (\cdot) $ is twice differentiable, and these points
tend themselves to zero. It is shown that a set of the limit
matrices coincides with the second-order subdifferential,
introduced by the author, of a positively homogeneous  function of
the second order at the point zero. The author's first and second
subdifferentials are used to calculate the first and second
codifferential of a codifferentiable function $f(\cdot) $. The
second hypodifferential and hyperdifferential of a function
$f(\cdot) $ are evaluated up to equivalence. Finally, the theorem
is proved that every Lipschitz quasidifferentiable function, whose
subdifferential and superdifferential are upper semicontinuous as
SVMs, is continuously codifferentiable. The proved theorems, that
give the rules for calculating subdifferentials and
codifferentials, are important for practical optimization.

\end{document}